%% file: InjEuclideanStiefel.tex
\documentclass[onefignum,onetabnum]{siamonline220329}

\input{InjEuclideanStiefel_shared}

\ifpdf
\hypersetup{
  pdftitle={The Euclidean Stiefel injectivity radius},
  pdfauthor={R. Zimmermann, and J. Stoye}
}
\fi

\begin{document}

\maketitle

% REQUIRED
\begin{abstract}
The injectivity radius of a manifold is an important quantity, both from a theoretical point of view and in terms of numerical applications.
It is the largest possible radius within which all geodesics are unique and length-minimizing.
In consequence, it is the largest possible radius within which calculations in Riemannian normal coordinates are well-defined.
A matrix manifold that arises frequently in a wide range of practical applications is the compact Stiefel manifold of orthogonal $p$-frames in $\R^n$. We observe that geodesics on this manifold are space curves of constant Frenet curvatures. Using this fact, we prove that the injectivity radius on the Stiefel manifold under the Euclidean metric is $\pi$.
\end{abstract}

% REQUIRED

% REQUIRED
\begin{MSCcodes}
  15B10, %Orthogonal matrices
  15B57, %Hermitian, skew-Hermitian, and related matrices
  65F99, % Numerical linear algebra, None of the above, but in this section
  53C30, %Differential geometry of homogeneous manifolds
  53C80, %Applications of global differential geometry to the sciences
\end{MSCcodes}

\section{Introduction}
The set of orthonormal $p$-frames in $\R^n$, i.e., the set of column-ortho-normal $(n\times p)$-matrices forms a Riemannian manifold, called the Stiefel manifold $St(n,p)$. Stiefel manifolds feature in a large variety of application problems, ranging from optimization
\cite{Absilbook:2008, boumal2023, sato2021} over numerical methods for differential equations \cite{BennerGugercinWillcox2015, Celledoni_2020, HueperRollingStiefel2008, Zim21}
to applications in statistics and data science \cite{Turaga_2008, Chakraborty2018, FLETCHER2020}. 
On a manifold, the selected Riemannian metric determines how length and angles are measured, and thus how geodesics are defined.
Geodesics are manifold curves that do not exhibit any intrinsic acceleration and are thus generalizations of straight lines in Euclidean spaces.
The injectivity radius is the largest possible radius within which geodesics are unique and length-minimizing. This is a strong geometric property. Among other things, it ensures that data processing operations in Riemannian normal coordinates are well-defined.
For general manifolds, it is very difficult to explicitly compute the injectivity radius. However, for a few special examples, such as the $n$-sphere or the Grassmann manifold of $p$-dimensional subspaces \cite{Wong1967, Wong1968}, the explicit number is known.

In this work, we add the Stiefel manifold to the list of manifolds with known injectivity radius. More precisely, we prove that the injectivity radius of the Stiefel manifold under the Euclidean metric is $\pi$,
\begin{equation*}
    \mathrm{inj}(St(n,p)) = \pi.
\end{equation*}
By a classical result from Riemannian geometry, the injectivity radius is related to the sectional curvatures and shortest closed geodesics on a manifold. In \cite{zimmermannstoye:2024}, sharp bounds on the sectional curvatures are provided. Here, we complete the picture by investigating shortest closed geodesics. The idea is to view the Euclidean Stiefel geodesics as curves in $\R^{np}$ and show that their Frenet curvatures are constant. Exploiting the normal form for constant-curvature curves of \cite[2.16 Remark]{kuhnel:2015}, we conclude that closed Euclidean Stiefel geodesics are at least of length $2\pi$ (\Cref{thm:shortestclosedgeos}). Since examples of explicit closed geodesic of length $2\pi$ exist, it follows that the injectivity radius of the Stiefel manifold under the Euclidean metric is $\pi$.

\paragraph{Organisation of the paper}
\Cref{sec:background} provides the required background on the Stiefel manifold, the injectivity radius and space curves with their Frenet curvatures. Our main results on closed geodesics and the injectivity radius of the Stiefel manifold are in \Cref{sec:main}. \Cref{sec:conclusions} concludes the paper.
\paragraph{Notation}
For the reader's convenience, we list the main acronyms and variables.

\begin{small}
\begin{tabular}{@{}ll@{}}
%\toprule
Symbol &  meaning \\
\midrule
$I,I_n$   & identity matrix, provided with a dimensional index if required \\
$\langle \cdot,\cdot\rangle$ & Euclidean inner product $\langle X, Y\rangle = \tr(X^T Y)$\\
$\|\cdot\|$ & (Frobenius) norm associated with the Euclidean inner product\\
$O(n)$    & orthogonal group $O(n) = \{Q\in \R^{n\times n}\mid Q^TQ = I_n\}$\\
$\Skew(n)$ & vector space of skew-symmetric matrices $\{A\in \R^{n\times n}\mid A^T = -A\}$.\\
$\mathcal{M}$ & a Riemannian manifold\\
$St(n,p)$ & Stiefel manifold $St(n,p)=\{U\in\R^{n\times p}\mid U^TU = I_p\}$\\
$T_USt(n,p)$ & tangent space of $St(n,p)$ at $U\in St(n,p)$\\
$\colspan(X)$ & span of a set of vectors or range of a matrix $X$\\
$\exp_m$ & matrix exponential $\exp_m(X) = \sum_{k=0}^\infty \frac{1}{k!}X^k$\\
$\kappa_j(t)$ & $j$th Frenet curvature of a regular space curve
\end{tabular}
\end{small}

\section{Background}
\label{sec:background}

We begin by introducing the Stiefel manifold and outline basic concepts of differentiable manifolds like geodesics, curvature and the injectivity radius in \Cref{sec:backgroundStiefel}. Related references are \cite{Zim21,gallier2011geometric}.
In \Cref{sec:backgroundCurves}, we provide a brief review of the essentials of curves in $\R^m$, focusing in particular on their Frenet curvatures. 

\subsection{The Stiefel manifold}
\label{sec:backgroundStiefel}

The Stiefel manifold is the compact, homogeneous matrix manifold of rectangular column-orthogonal matrices 
\begin{equation*}
St(n,p) := \{U\in\R^{n\times p}|U^TU = I_p\},
\end{equation*}
for $n\geq p$. It is an $np-\frac12(p(p+1))$-dimensional embedded submanifold of $\R^{np}\cong\R^{n\times p}$. 
The Stiefel manifold is extensively studied in the literature, see e.g., 
\cite{EdelmanAriasSmith:1999, Absilbook:2008, Zim21}.
For any point $U\in St(n,p)$ on the Stiefel manifold, the tangent space of $St(n,p)$ at $U$ is given by 
\begin{equation*}
T_U St(n,p) = \{\Delta\in\R^{n\times p}| U^T\Delta = -\Delta^TU\}\subset \R^{n\times p}.
\end{equation*}
Stiefel tangent vectors are rectangular matrices $\Delta = \mathbf{U}\begin{pmatrix}
    A\\
    B
\end{pmatrix} = UA+U_\bot B\in\R^{n\times p}$, where $\mathbf{U}=\begin{pmatrix}
    U& U_\bot
\end{pmatrix}\in O(n)$, $A\in \Skew(p), B\in \R^{(n-p)\times p}$.  A Riemannian metric on the Stiefel manifold is obtained from restricting the Euclidean inner product of $\R^{n\times p}$ to the Stiefel tangent spaces. This yields the so-called {\em Euclidean metric} on $St(n,p)\subset\R^{n\times p}$,
\begin{equation}\label{eq:EuclideanMetric}
\langle \Delta,\tilde \Delta\rangle^{St}_e
= \tr\left(\begin{pmatrix}
    A^T & B^T
\end{pmatrix}\mathbf{U}^T \mathbf{U}
\begin{pmatrix}
    \tilde A\\
    \tilde B
\end{pmatrix}\right)
= \tr(A^T\tilde A) + \tr(B^T\tilde B).
\end{equation}
%Verweis auch zu family of metric?
The Euclidean metric is also included in a parametric family of metrics considered in \cite[Def. 3.1, Section 4.2]{HueperMarkinaLeite2020}.

The length of a tangent vector $\Delta\in T_USt(n,p)$ is $\|\Delta\|_U = \sqrt{\langle \Delta,\Delta\rangle_U}$ and the length of a curve $c\colon \left[a,b\right]\to St(n,p)$ is
\begin{equation*}
L(c) := \int_a^b \Vert\dot{c}(t)\Vert_{c(t)}dt.
\end{equation*}
In the following, we often omit the index that indicates the dependence of the metric on the base point. Moreover, $\|\cdot\|$ always denotes the Euclidean metric (with associated Frobenius norm), regardless if  applied on $\R^{n\times p}$ or on a tangent space $T_USt(n,p).$
The Riemannian distance between two points $U,\,\tilde U\in St(n,p)$ is 
\begin{equation*}
\text{dist}_{St(n,p)}(U,\tilde U) := \inf\{L(c\vert_{[a,b]})\mid c(a)=U, c(b) = \tilde U\}.
\end{equation*}
Geodesics are candidates for length-minimizing curves. They are charactarized by the fact that they have no intrinsic, or, more precisely, no covariant acceleration and can be obtained as local solutions to an ordinary differiential equation, see \cite[5.17 Definition, 5.18 Corollary]{kuhnel:2015}. The Riemannian exponential map on a Riemannian manifold $\mathcal{M}$ is based on geodesics. For a point $p\in\mathcal{M}$, the Riemannian exponential sends a tangent vector $v\in T_p\mathcal{M}$ to the endpoint of the geodesic on the unit interval that starts from $p$ with velocity $v$.
For a precise definition, see \cite[5.19 Definition]{kuhnel:2015} or \cite[p.72]{Lee:RiemannianManifolds} or any other textbook on Riemannian geometry.

%let $\mathcal{M}$ be a Riemannian manifold.
%Let $\gamma_{p,v}$ be the geodesic starting from $p\in \mathcal{M}$ with velocity $v\in T_p\mathcal{M}$. The Riemannian exponential is defined as
%\begin{equation*}
%\text{Exp}_p\colon T_p \mathcal{M}\supset B_{\epsilon}(0)\to %\mathcal{M},\quad v\mapsto q := \gamma_{p,v}(1).
%\end{equation*}
%For technical reasons, $\epsilon>0$ must be small enough so that $\gamma_{p,v}(t)$ is defined on the unit interval $\left[0,1\right]$.
%So we can go to 1 and not experience some ununiqueness before

%
A formula for the Riemannian exponential of the Stiefel manifold under the Euclidean metric \eqref{eq:EuclideanMetric} is given in \cite[Prop. 1 ($\alpha = -\frac12$)]{ZimmermannHueper2022}:
    For $U\in St(n,p)$ and $\Delta\in T_U St(n,p)$, it  reads
    \begin{equation}\label{lem:ExplicitExponential}
        \text{Exp}_U(\Delta) =  \begin{pmatrix}
        U & Q
    \end{pmatrix}
    \exp_m\left(\begin{pmatrix}
        2A & -B^T\\
        B  & \mathbf{0} 
    \end{pmatrix}\right)
    \begin{pmatrix}
      I_p\\
      \mathbf{0}
    \end{pmatrix}\exp_m(-A)\in \R^{n\times p},
    \end{equation}
    with $A = U^T\Delta\in\text{skew}(p)$ and $QB = (I-UU^T)\Delta\in\R^{n\times p}$ being any matrix decomposition with $Q\in St(n,p)$ and $B\in\R^{p\times p}$.
The corresponding geodesic is $\gamma: t\mapsto \text{Exp}_U(t\Delta)$, and the length of the geodesic on $[0,1]$ is $L(\gamma\big\vert_{[0,1]}) = \|\Delta\|$.

There is a clear criterion under which a geodesic is length-minimizing. On a Riemannian manifold $\mathcal{M}$, the injectivity radius is the largest possible radius within which geodesics are unique and length-minimizing, regardless of where you start from.
In loose words, as long as you stay within the injectivity radius when travelling along a geodesic, you are guaranteed not to travel unnecessary distances.
\begin{definition}[Injectivity radius]\label{def:injrad}
Let $\epsilon$ be the maximum radius of $B_{\epsilon}(0)$ such that the Riemannian exponential at $p$, $\text{Exp}_p\colon T_p \mathcal{M}\supset B_{\epsilon}(0)\to D_p\subset\mathcal{M}$, is invertible. 
Then, $\epsilon$ is called the injectivity radius of $\mathcal{M}$ at $p$ and is denoted by $\mathrm{inj}_p(\mathcal{M})$.\\
The infimum of $\mathrm{inj}_p(M)$ over all $p\in\mathcal{M}$ is called the (global) injectivity radius of $\mathcal{M}$,
\[\mathrm{inj}(\mathcal{M}) = \inf_{p\in\mathcal{M}}\mathrm{inj}_p(\mathcal{M}).\]
\end{definition}
As with the Riemannian exponential, the injectivity radius depends on the underlying metric of the manifold. For details, we refer the reader to \cite[Chap. 13]{DoCarmo2013riemannian}.

The injectivity radius is closely related to the sectional curvatures of the manifold, as we will see later. For two linearly independent tangent vectors at a point, one can consider their span as a subplane of the tangent space. The manifold image of this coordinate plane under the Riemannian exponential is then a two-dimensional submanifold, and the associated sectional curvature quantifies the deviation of this submanifold from being flat.

 In \cite{nguyen2022curvature}, the sectional curvatures of the Stiefel manifold are computed for a family of Riemannian metrics from \cite{HueperMarkinaLeite2020}.
In \cite{zimmermannstoye:2024}, sharp bounds on the sectional curvature for almost all dimensions of the Stiefel manifold are proven. The essential finding is in the next theorem.
\begin{theorem}{\cite[Theorem 10]{zimmermannstoye:2024}}
\label{thm:StiefelCurvBound_eucl}
	The sectional curvatures $K^{St}_e$ on the Stiefel manifold $St(n,p)$, $n\geq p$ under the Euclidean metric are globally bounded by 
	\[
	-\frac 12 \leq K^{St}_e \leq 1.
	\]
 \iffalse 
 \begin{enumerate}
     \item For $p\geq 2, n\geq p+2$, the bound is sharp.
    \item For $p=1$ and $n\geq 3$, the sectional curvature has a constant value of 
    $$K^{St}_e\equiv 1.$$ 
    \item For $p=n\geq 4$, it holds
    	\[
	0 \leq K^{St}_e(X,Y) \leq \frac14.
	\]
       \item For $p=n=3$, it holds
    	\[
	K^{St}_e(X,Y) \equiv \frac18.
	\]
    \item For $(n,p)=(3,2)$ the sectional curvature is bounded by
    \[ 
        -\frac12\leq K^{St}_e(X,Y)\leq \frac12.
    \]
    The bounds are sharp.
    The upper bound is attained, e.g., for
    \[
    X = \begin{pmatrix}
        0&0\\
        0&0\\
        -1&0
    \end{pmatrix}, \quad 
    Y = \begin{pmatrix}
        0&-\frac{1}{\sqrt{2}}\\
        \frac{1}{\sqrt{2}}&0\\
        0&0
    \end{pmatrix}.
    \]
        \item For $n\geq 4$ and $(n,p)=(n,n-1)$ the sectional curvature is bounded by
    \[ 
        -\frac12\leq K^{St}_e(X,Y)\leq \frac23.
    \]
    The lower bound is sharp.
 \end{enumerate}
 \fi
 \end{theorem}
For certain dimensions, the upper bound might be smaller (but always strictly positive). For example, $K^{St}_e(X,Y) \equiv \frac18$ on $St(3,3)$.
This is detailed in \cite[Theorem 10]{zimmermannstoye:2024}. Yet, these refinements have no impact on the considerations here. The reason for this is the Klingenberg theorem from Riemannian geometry, which relates the injectivity radius to the sectional curvatures and the length of shortest closed geodesics.
Eventually, it turns out that for a Stiefel manifold with any upper curvature bound $0<C\leq 1$, it is always the closed geodesics that determine the injectivity radius.  
\begin{theorem}[Klingenberg, stated as Lemma 6.4.7 in \cite{petersen2016riemannian}]
\label{thm:Klingenberg}
	Let $\mathcal{M}$ be a compact Riemannian manifold with sectional curvatures bounded by $C>0$. Then the injectivity radius $\mathrm{inj}(p)$ at any $p\in \mathcal{M}$ satisfies
	\[
	\mathrm{inj}(p)\geq \min\left\{\frac{\pi}{\sqrt{C}}, \frac12 l_p\right\},
	\]
	where $l_p$ is the length of a shortest closed geodesic starting from $p$.
	For the global injectivity radius, it holds
		\[
	\mathrm{inj}(\mathcal{M})  \geq\frac{\pi}{\sqrt{C}}
	\quad \text{ or } \quad \mathrm{inj}(\mathcal{M}) =\frac12 l,
	\]
	where $l$ is the length of a shortest closed geodesic on $\mathcal{M}$.
\end{theorem}
% 
%For the Stiefel manifold, the curvature condition of this theorem is knwon from \Cref{thm:StiefelCurvBound_eucl}. In \Cref{sec:main}, we will concentrate on finding the length of shortest closed geodesics on the Stiefel manifold under the Eucliden metric. 
%
%
%
%
%
\subsection{Space curves and their Frenet curvatures}\label{sec:backgroundCurves}
We provide a short review of the essentials on curves in $\R^m$ taken from the textbooks \cite[Sections 1.2, 1.3]{Klingenberg:1978} and \cite[Chapter 2]{kuhnel:2015}.
Related references are \cite{Gluck:1966, Monterde:2005}.

For any space curve $c:I\to \R^N$ with linear independent derivatives $\dot c(t), \ddot c(t), \ldots, c^{(N-1)}(t)$, there exists a unique orthonormal moving $N$-frame $e_1(t), \ldots, e_N(t)$, called the {\em distinguished Frenet frame associated with} $c$ such that:
\begin{enumerate}
    \item Each $t\mapsto e_j(t)\in\R^N$ is a vector field along $c$.
    \item $\langle e_j(t),e_k(t)\rangle=\delta_{jk}$.
    \item For $k=1,\ldots,N-1$,   $c^{(k)}(t)\in 
\colspan\{e_1(t), \ldots, e_k(t)\}$.
    \item For $k=1,\ldots,N-1$, the sets of vectors
    $\dot c(t), \ddot c(t), \ldots, c^{(k)}(t)$
    and $e_1(t), \ldots, e_k(t)$ have the same orientation.
    \item The frame $\left(e_1(t), \ldots, e_N(t)\right)$ is of positive orientation.
\end{enumerate}
(cf. \cite[Def. 1.2.1, Prop. 1.2.2]{Klingenberg:1978}.)
If there is $k<N-1$ such that $c^{(k)}(t)$ is in the span of $\dot c(t), \ddot c(t), \ldots, c^{(k-1)}(t)$ on a small interval,
then this also holds for all higher derivatives. Hence, the curve is completely contained in a $k$-dimensional subspace of $\R^N$ and the analysis can be reduced to this subspace. 
For presenting the basic theory, it is no loss of generality to assume that $\dim(\colspan(\dot c(t), \ddot c(t), \ldots, c^{(N-1)}(t))) = N-1$ on a suitable interval of definition.

If the curve $c$ is parameterized by its arc length, we have $\|\dot c(t)\|\equiv 1$ and $e_1(t) = \dot c(t)$.
The coordinates of the distinguished Frenet frame for an arc-length curve $c$ satisfy a highly structured ODE:
\begin{equation}
\label{eq:Frenet-ODE}
\begin{pmatrix}
    e_1(t)\\
    e_2(t)\\
    \vdots\\
    e_{N-1}(t)\\
    e_{N}(t)
\end{pmatrix}'
=
\begin{pmatrix}
        0    & \kappa_1(t) &       &                &  \\
-\kappa_1(t) & 0           &\ddots &                &   \\
             &\ddots       &\ddots &\ddots          &   \\
             &             &\ddots & 0              &  -\kappa_{N-1}(t)\\
             &             &       &-\kappa_{N-1}(t)& 0
\end{pmatrix}
\begin{pmatrix}
    e_1(t)\\
    e_2(t)\\
    \vdots\\
    e_{N-1}(t)\\
    e_{N}(t)
\end{pmatrix}.
\end{equation}
The $j$th coefficient function is 
\[
    \kappa_j(t) = \langle \dot e_j(t), e_{j+1}(t)\rangle.
\]
Adopting the terminology from \cite{kuhnel:2015}, $t\mapsto \kappa_j(t)$ is called the $j$th {\em Frenet  curvature} of $c$. An alternative name is `$j$th Euclidean curvature'.
The curvature coefficients determine the curve $c$ uniquely up to isometries -- via the ODE system \eqref{eq:Frenet-ODE}, \cite[Theorems 1.3.5 \& 1.3.6]{Klingenberg:1978}.

If all curvature coefficients are constant, the system \eqref{eq:Frenet-ODE} can be integrated in closed form. The solution curve has the shape of a twisted geodesic on a flat torus, see \cite[2.16 Remark]{kuhnel:2015} or \cite[Corollary 1]{Monterde:2005}.
This means that in even dimensions $N=2m$, $c$ features the expression
\begin{equation}
    \label{eq:curvconstcurv_even}
    c(t) = \left(a_1\cos(b_1t),a_1\sin(b_1t), \ldots, a_m\cos(b_mt),a_m\sin(b_mt)\right)^T\in \R^{N}.
\end{equation}
If $N$ is odd, $N = 2m+1$,
\begin{equation}
    \label{eq:curvconstcurv_odd}
    c(t) = \left(a_1\cos(b_1t),a_1\sin(b_1t), \ldots, a_m\cos(b_mt),a_m\sin(b_mt), a_{m+1} t\right)^T\in \R^{N}.
\end{equation}
For a true Frenet curve (with linearly independent derivatives), the coefficients $b_i$ are pairwise distinct.
Otherwise, the curve can be reduced to a lower-dimensional subspace to become a true Frenet curve therein.
\section{Main results}
\label{sec:main}
In this section, we prove that the injectivity radius of the Stiefel manifold under the Euclidean metric is $\pi$. The line of argumentation is as follows:
First, we show that up to a fixed coordinate frame, geodesics on $St(n,p)\subset \R^{n\times p}$ are curves of constant Frenet curvatures in $\R^{2p\times p}\cong \R^{2p^2}$ and thus take the form of \eqref{eq:curvconstcurv_even}. Using curvature bounds, we show that under the Euclidean metric, the shortest closed curves of this type have a length of $2\pi$.
\begin{lemma}
\label{lem:ConstCurvSpaceCurves}
Let $\gamma:I\to \R^N$ be a smooth space curve such that for all $j\in\N$, the $j$th derivative of $\gamma$ is of constant norm, $\|\gamma^{(j)}(t)\|\equiv \mathrm{const}$.
Then $\gamma$ is a curve of constant Frenet curvatures.
\end{lemma}
\begin{proof}
We start by showing that 
  all inner products $\langle \gamma^{(j)}(t),\gamma^{(k)}(t)\rangle$ are constant in time. This can be established by induction.
The start is obvious, so suppose that the claim has been checked for all inner products of derivatives up to order $k$.
Consider $\langle \gamma^{(k+1)}(t),\gamma^{(j)}(t)\rangle$ and distinguish two cases. First, let $j<k$. Then
\[
    \langle \gamma^{(k+1)}(t),\gamma^{(j)}(t)\rangle   
    = \frac{d}{dt} \langle \gamma^{(k)}(t),\gamma^{(j)}(t)\rangle
    -\langle \gamma^{(k)}(t),\gamma^{(j+1)}(t)\rangle = 0 - \text{const}.,
\]
because both terms on the right are constant by the induction hypothesis.
For $j=k$, 
\[
    \langle \gamma^{(k+1)}(t),\gamma^{(k)}(t)\rangle 
    = \frac12 \frac{d}{dt} \langle \gamma^{(k)}(t),\gamma^{(k)}(t)\rangle =
    \frac{d}{dt} \text{const} = 0.
\]
In fact, much more detailed information\footnote{e.g., $\langle \gamma^{(2k+1)}, \gamma\rangle = 0, \langle \gamma^{(2k)}, \gamma^{(1)}\rangle = 0, \ldots$} about the inner products can be obtained, but this is not required for the argument that is to follow.
The distinguished Frenet frame is obtained from a Gram-Schmidt process applied to $\dot\gamma, \ddot\gamma,\ldots,\gamma^{(N-1)}$. 
As mentioned in \Cref{sec:backgroundCurves}, We can restrict the considerations to the case where the derivative vectors are pointwise linearly independent. 
%\TODO{If linearly indep. at $t_0$, then also lin indep on $[t_0-\epsilon, t_0+\epsilon]$; argue locally as with the solutions to ODEs. Yet no problem on Stiefel: all derivatives are either linear independent for all $t$, or linear dependent for all $t$, if the constant matrices $[M_j^T,N_j^T]^T$ from \eqref{eq:geo_jth_deriv} are.} 
If the columns of $E_{i-1}(t) = \left(e_1(t),\ldots,e_{i-1}(t)\right)$ form already an orthonormal basis of 
$\text{span}(\dot\gamma(t), \ddot\gamma(t),\ldots,\gamma^{(i-1)}(t))$, the $i$th step in the Gram-Schmidt process reads
\begin{eqnarray*}
    \tilde e_i(t) &=& \left(I-E_{i-1}(t)E_{i-1}^T(t)\right) \gamma^{(i)}(t) 
    = \gamma^{(i)}(t) - E_{i-1}(t)
    \begin{pmatrix}
    \langle e_1(t), \gamma^{(i)}(t)\rangle    \\
    \vdots\\
    \langle e_{i-1}(t), \gamma^{(i)}(t)\rangle
    \end{pmatrix}, \quad
     e_i(t) = \frac{\tilde e_i(t)}{\| \tilde e_i(t) \|}.
\end{eqnarray*}
Because all inner products between the derivatives are constant, all coefficients in the Gram-Schmidt process are constant and the orthonormalization is eventually realized by an upper triangular matrix $R$ with constant coefficients
\[
     \left(\dot\gamma(t), \ddot\gamma(t),\ldots,\gamma^{(N-1)}(t)\right)=: \Gamma(t) 
    = E_{N-1}(t) R.
\]
As a consequence,
\[
    \dot E_{N-1}(t)^T E_{N-1}(t) = R^{-T} \dot \Gamma(t)^T \Gamma(t) R^{-1}  = \text{Const.}
\]
is a constant matrix. This is because $\dot \Gamma(t)^T \Gamma(t)$ only contains inner products of derivatives of $\gamma$, which are constant. In particular, all the Frenet curvatures 
$\kappa_j(t) = \langle \dot e_j(t), e_{j+1}(t)\rangle \equiv \kappa_j$
are constant for $j=1,\ldots, N-2$.\\
For $j=N-1$, we have to show that $\langle \dot e_{N-1}(t), e_N(t) \rangle \equiv \text{const.}$. Here, $e_N(t)$ is the pointwise unique vector that has been added to $\text{span}\left(\dot\gamma(t),\ldots,\gamma^{(N-1)}(t)\right) = \text{span}\left(e_1(t),\ldots,e_{(N-1)}(t)\right)$ so that the Frenet frame $E_N(t)$ has determinant one (positive orientation). If there is $t_0$ such that $\gamma^{(N)}(t_0)$ is linearly independent from the other derivatives $\gamma^{(j)}(t_0),\,j=1,\ldots,N-1$,
then this is also the case on a small interval around $t_0$.
Then the $t$-dependent Gram-Schmidt process can be completed with the last vector function $e_N(t)$ given as a linear combination of the $\gamma^{j}(t),\,j=1,\ldots,N$ and the proper choice of the sign. Then $\kappa_{N-1}(t) = \langle \dot e_{N-1}(t), e_{N}(t)\rangle$
is constant by the same argument as above. 

If there is no such $t_0$, then $\gamma^{(N)}(t)$ is linearly dependent on $\gamma^{(j)}(t), j=1,\ldots,N-1$ on a suitable interval. In this case, $\dot e_{N-1}(t)= \dot \Gamma(t) \begin{pmatrix}
    r_{1,N-1}\\
    \vdots\\
    r_{N-1,N-1}
\end{pmatrix}\in \text{span}\{\gamma^{(j)}(t)\mid j=2,\ldots,N\}$ is orthogonal to $e_N(t)$.
Hence $\kappa_{N-1}(t)\equiv 0$.
\end{proof}
\begin{lemma}
\label{lem:ConstCurvStiefel}
Euclidean Stiefel geodesics are curves of constant Frenet curvatures.
\iffalse
\RZ{Eigentlich gilt das genauso für alle $\alpha$-Metriken aus \cite{HueperMarkinaLeite2020}, also auch die kanonische. Ich war schon dabei, das Resultat umzuschreiben, habe dann aber doch gezögert.}\JS{Ja, das habe ich auch schon gedacht. Das sind ja auch geschlossene Kurven im $R^{np}$ und die Formel für die Geodäten unterscheiden sich nur um einen von $\alpha$ abhängigen Faktor vor dem $A$. Ich frage mich, ob wir nicht damit eine Aussage über die Länge einer geschlossener Geodäten bzgl einer beliebigen $\alpha$-Metriken bekommen. (Mit der Beispielgeodäten haben wir wegen $A = 0$ für jedes $alpha$ eine geschlossene Geodäte der Länge $2\pi$.)}
\RZ{Man muss nur aufpassen, wie man dann Längen misst. Die Normalform mit (sin, cos)-Komponenten wird ja durch Euklidische Isometrien erreicht.}
\JS{Ja gut, das könnte Probleme geben, das stimmt...}
\fi
%The Euclidean metric and the canonical metric are included in this family of metrics and correspond to $\alpha = -\frac12$ (Euclidean case) and $\alpha =0$ (canonical case).
\end{lemma}
\begin{proof}
   We start from the matrix formula for the Euclidean Stiefel geodesics restated in  \eqref{lem:ExplicitExponential}.
   Given $U\in St(n,p)$ and $\Delta\in T_USt(n,p)$, compute
   $A =U^T\Delta\in \Skew(p)$, and a decomposition
   $QB = (I-UU^T)\Delta$ with $Q\in St(n,p), B\in \R^{p\times p}$. Then the associated geodesic is
\[
    \gamma(t) = 
    \begin{pmatrix}
        U & Q
    \end{pmatrix}
    \exp_m\left(t\begin{pmatrix}
        2A & -B^T\\
        B  & \mathbf{0} 
    \end{pmatrix}\right)
    \begin{pmatrix}
      I_p\\
      \mathbf{0}
    \end{pmatrix}\exp_m\left(-t A\right)\in \R^{n\times p}.
\] 
Since $\begin{pmatrix}
        U & Q
    \end{pmatrix}$
is nothing but a fixed coordinate frame,
it has no effect on the following considerations and we are free to drop it and consider $\gamma(t)$ as a space curve in $\R^N$, $N=2p^2$, whenever useful.
The $j$th derivative of $\gamma$ is then
\begin{equation}
 \label{eq:geocurve_jth_deriv}
    \gamma^{(j)}(t) = 
    \exp_m\left(t \begin{pmatrix}
        2A & -B^T\\
        B  & \mathbf{0} 
    \end{pmatrix}\right)
    \begin{pmatrix}
      M_j\\
      N_j
    \end{pmatrix}\exp_m(-tA), \quad j=0,1,2,\ldots
\end{equation}
where $\begin{pmatrix}
      M_j\\
      N_j
    \end{pmatrix}$ is given iteratively by
\begin{equation}
  \begin{pmatrix}
      M_0\\
      N_0
  \end{pmatrix}  =
\begin{pmatrix}
      I_p\\
      \mathbf{0}
    \end{pmatrix}, \quad
      \begin{pmatrix}
      M_{j+1}\\
      N_{j+1}
    \end{pmatrix}
    =
    \begin{pmatrix}
        2A & -B^T\\
        B  & \mathbf{0} 
    \end{pmatrix}
          \begin{pmatrix}
      M_j\\
      N_j
    \end{pmatrix}
    -
              \begin{pmatrix}
      M_jA\\
      N_jA
    \end{pmatrix} \quad j=0,1,2,\ldots .
\label{eq:geo_jth_deriv}
\end{equation}
The important thing to note is that for all $j$,
\begin{equation}
\label{eq:const_deriv_norm}
    \|\gamma^{(j)}(t)\|^2 = \|\begin{pmatrix}
      M_j\\
      N_j
    \end{pmatrix}\|^2 = \tr(M_j^TM_j + N_j^TN_j) = \text{const.}
\end{equation}
The claim follows from \Cref{lem:ConstCurvSpaceCurves}.
\end{proof}
\begin{remark}
    Thanks to the skew-symmetric structure of the matrices in the geodesic formula of \cite[Prop. 1]{ZimmermannHueper2022}, the result was rather obvious. It is much less clear when looking at the original formula for Euclidean Stiefel geodesics from \cite[p. 310]{EdelmanAriasSmith:1999}.\\
    Essentially, the same argument holds for the Stiefel geodesics under any metric from the one-parameter family of $\alpha$--metrics from \cite{HueperMarkinaLeite2020}. 
    The only difference is that in the explicit form of the matrix blocks $M_j,N_j$, a constant factor associated with the metric's parameter $\alpha$ appears. Hence, when considered as space curves in $\R^{2p^2}$, all Stiefel geodesics are curves of constant Frenet curvatures, regardless of the chosen $\alpha$--metric.
\end{remark}
\begin{corollary}
\label{cor:normal_form_const_curve}
Stiefel geodesics under the Euclidean metric have the normal form
\begin{equation}
    \label{eq:StiefelGeo_even}
    \gamma(t) = \left(a_1\cos(b_1t),a_1\sin(b_1t), \ldots, a_m\cos(b_mt),a_m\sin(b_mt)\right)^T\in \R^{2p^2}.
\end{equation}
\end{corollary}
 \begin{proof}
    Note that the dimension $N=2p^2$ is always even. The claim follows from the characterization of constant-curvature curves of \cite[Corollary 1]{Monterde:2005} that was stated here as equation \eqref{eq:curvconstcurv_even}.
 \end{proof}

\begin{theorem}\label{thm:shortestclosedgeos}
  The shortest closed geodesics on $St(n,p)$ under the Euclidean metric have a length of $2\pi$.
\end{theorem}
\begin{proof}
    Let $\gamma$ be a closed geodesic on $St(n,p)$ under the Euclidean metric.
    W.l.o.g. we assume that $\gamma$ is paramaterized by the arc length, i.e., $\|\dot\gamma(t)\|\equiv 1$. Moreover, as argued in \Cref{lem:ConstCurvStiefel}, it is sufficient to consider $\gamma$ as a curve in $\R^{2p\times p} \cong \R^{N}$, $N=2p^2$.
    Using the form of \eqref{eq:geocurve_jth_deriv} for $\gamma$, we have
    \begin{eqnarray*}
        \|\dot\gamma(t)\|^2 &=& \tr(A^TA) + \tr(B^TB) = 1,\\
        \|\ddot\gamma(t)\|^2 &=& \kappa_1^2 = \|-A^TA-B^TB\|_2^2 =  \tr(A^4) + \tr((2A^TA +B^TB)B^TB) \\
        &\leq& \|A\|^4 + \|2A^TA + B^TB\|\|B^TB\| \leq  \|A\|^4 + (\|A\|^2 + \underbrace{\|A\|^2 +  \|B\|^2}_{=1})\|B\|^2\\
        &=& \|A\|^4 + (\|A\|^2 +1)(1-\|A\|^2) = 1.
    \end{eqnarray*}
    Hence, the maximum possible curvature (in the sense of space curves) of a Stiefel geodesic is one.
    Now, transform $\gamma$ to the normal form of \eqref{cor:normal_form_const_curve} via a Euclidean isometry. Then,
    \[ 
      \gamma(t) = \left(a_1\cos(b_1t),a_1\sin(b_1t), \ldots, a_m\cos(b_mt),a_m\sin(b_mt)\right)^T\in \R^{2p^2}, \quad m=p^2.
    \]
    In these coordinates, one sees that $\gamma$ can only close\footnote{Mind the analogy to Wong's characterization of closed geodesics on the Grassmann manifold in \cite{Wong1967}.}, if 
    there is a natural number $k\in \N$ such that for each $i=1,\ldots,m$,
    \[
        b_ik\in 2\pi\N, \quad \text{say } k=\frac{2\pi l_i}{b_i}, l_i\in \N.
    \]

    In this case, the length of the loop 
    $\gamma(0) = (a_1, 0, a_2,0,\ldots, a_m, 0)^T=\gamma(k)$ is
    $L(\gamma_{[0,k]}) = \int_0^k\|\dot\gamma(t)\|dt = k$.
    We have $1 = \|\dot\gamma(t)\|^2 = \sum_{i=1}^m a_i^2 b_i^2$ and
    $1\geq \kappa_1^2 = \|\ddot\gamma(t)\|^2 =\sum_{i=1}^m a_i^2 b_i^4$.
    If $|b_i|> 1$ for all $i$, then $1\geq \sum_{i=1}^m a_i^2 b_i^4> \sum_{i=1}^m a_i^2 b_i^2=1$, a contradiction.
    Hence, there is $b_i\leq 1$, which gives
    \[
        L(\gamma|_{[0,k]}) = k = \frac{2\pi l_i}{b_i} \geq 2\pi.
    \]
\end{proof}

A closed unit-speed geodesic of shortest possible length $2\pi$ on the Stiefel manifold $St(4,2)$ under the Euclidean metric is
\[
    \gamma(t) = \exp_m\left(t\begin{pmatrix}
        2A & -B^T\\
        B  & \mathbf{0} 
    \end{pmatrix}\right)
    \begin{pmatrix}
      \exp_m(-tA)\\
      \mathbf{0}
    \end{pmatrix}, \quad A = \mathbf{0}\in \Skew(2), 
    \quad B = \begin{pmatrix}
        1 & 0\\
        0 & 0
    \end{pmatrix}.
\]
For this specific choice of matrices $A,B$, the matrix exponential evaluates to 
$$\exp_m\left(t\begin{pmatrix}
        \mathbf{0} & -tB^T\\
        tB  & \mathbf{0} 
    \end{pmatrix}\right) 
        \begin{pmatrix}
      I_2\\
      \mathbf{0}
    \end{pmatrix}
    = \begin{pmatrix}
        \cos(tB) & -\sin(tB)\\
        \sin(tB)  & \cos(tB)
    \end{pmatrix}\begin{pmatrix}
      I_2\\
      \mathbf{0}
    \end{pmatrix}=
    \begin{pmatrix}
      \cos(tB)\\
      \sin(tB)
    \end{pmatrix}
    =\begin{pmatrix}
        \cos(t) & 0\\
        0       & 1\\
        \hline
        \sin(t) & 0\\
         0      & 0
    \end{pmatrix}
    .$$
The loop starts and closes at $\gamma(0) = \gamma(2\pi) = 
\begin{pmatrix}
I_2\\
\mathbf{0}
\end{pmatrix}$. 
Identifying $\R^{2p\times p} \cong \R^8$ and rearranging the order of the coordinates,
this loop can be written as $\gamma(t) = (\cos(1\cdot t),\sin(1\cdot t), 0,0, 1,0, 0,0)^T$, 
which is simply a planar unit circle in a certain coordinate plane.
Of course, this geodesic can be embedded in higher dimensions and can easily be adapted for all other combinations of $n$ and $p$.
More examples of closed geodesics and their length with respect to the metrics of the family from 
\cite{HueperMarkinaLeite2020} are given in \cite[Section 6]{absil2024ultimate}.

\begin{theorem}\label{thm:injEuclidStiefel}
  The injectivity radius of the  Stiefel manifold $St(n,p)$ under the Euclidean metric is $\pi$. 
\end{theorem}
We emphasize that the theorem holds for all dimensions, including the special cases of $(n,p)=(n,1)$,  $(n,p)=(n,n)$ and $(n,p) = (n,n-1)$.
\begin{proof}
The sectional curvature of $St(n,p)$ is bounded from above by one by \Cref{thm:StiefelCurvBound_eucl}.
The length of a shortest closed geodesic is at least $l=2\pi$ by \Cref{thm:shortestclosedgeos}. Examples of geodesics of length $2\pi$
(which are essentially planar circles, possibly embedded in higher dimensions) exist on all Stiefel manifolds.
    The theorem follows from combining these facts with \Cref{thm:Klingenberg}.
\end{proof}

\section{Summary and outlook}
\label{sec:conclusions}
The injectivity radius of the Stiefel manifold under the Euclidean metric equals $\pi$.
An elementary but essential observation is that all Stiefel geodesics are curves of constant Frenet curvatures
when considered as space curves in the embedding space.
The Euclidean framework makes it possible to transfer the geodesics directly to a certain Frenet normal form without distorting the curve's geometric characteristics.
Euclidean Stiefel geodesics have at most a (first Frenet) curvature of $\kappa_1 = 1$. Hence, they are at least as long as the planar circle of unit curvature, i.e., $2\pi$.

Because sharp bounds on the sectional curvature were known beforehand, this was the last building block to fully leveraging the classical Klingenberg's Theorem to conclude that the injectivity radius of the Stiefel manifold under the Euclidean metric is $\pi$ across all dimensions.

As future work, it is planned to transfer the investigation to closed Stiefel geodesics under the one-parameter family of Riemannian metrics from \cite{HueperMarkinaLeite2020}. This may complement the work of \cite{absil2024ultimate}.
Shortest closed geodesics are rather special objects of low intrinsic dimension. We plan to add numerical experiments to the theory to investigate the behavior of the Stiefel exponential along generic directions.

%Consequently, geodesics with lengths of at most $\pi$ are guaranteed to be unique and length-minimizing, and calculations within Riemannian normal coordinates remain well-defined within a radius of $\pi$. 

%\appendix
%\section{An example appendix} 

%\bibliographystyle{siamplain}
%\bibliography{references}

\end{document}

%% file: InjEuclideanStiefel_shared.tex
\usepackage{amsfonts}
\usepackage{graphicx}
\usepackage{epstopdf}
\usepackage{booktabs}
\usepackage{algorithmic}
\ifpdf
  \DeclareGraphicsExtensions{.eps,.pdf,.png,.jpg}
\else
  \DeclareGraphicsExtensions{.eps}
\fi

\usepackage{enumitem}
\setlist[enumerate]{leftmargin=.5in}
\setlist[itemize]{leftmargin=.5in}

\newsiamremark{remark}{Remark}
\newsiamremark{hypothesis}{Hypothesis}
\crefname{hypothesis}{Hypothesis}{Hypotheses}
\newsiamthm{claim}{Claim}

\newcommand{\N}{\mathbb N}

\newcommand{\R}{\mathbb R}

\newcommand{\colspan}{\operatorname{span}}

\newcommand{\tr}{\operatorname{tr}}

\newcommand{\Skew}{\operatorname{skew}} %\skew exists already in LaTeX

\newcommand{\diff}[1][t]{\mathrm{d}}
\newcommand{\Diff}[1][t]{\mathrm{D}}

\makeatletter
\renewcommand*\env@matrix[1][*\c@MaxMatrixCols c]{%
	\hskip -\arraycolsep
	\let\@ifnextchar\new@ifnextchar
	\array{#1}}
\makeatother

\newcount\Comments  % 0 suppresses notes to selves in text
\newcommand{\kibitz}[2]{\ifnum\Comments=1\textcolor{#1}{#2}\fi}

\newcommand{\RZ}[1]  {\kibitz{teal}   {[RZ: #1]}}
\newcommand{\JS}[1]  {\kibitz{blue}   {[JS: #1]}}

\headers{The Euclidean Stiefel injectivity radius}{R. Zimmermann and J. Stoye}

\title{The injectivity radius of the compact Stiefel manifold under the Euclidean metric%\thanks{Submitted to the editors DATE}
}

\author{Ralf Zimmermann\thanks{Department of Mathematics and Computer Science, University of Southern Denmark, Odense 
  (\email{zimmermann@imada.sdu.dk})}
\and Jakob Stoye\thanks{Institute for Numerical Analysis, TU Braunschweig, Germany 
  (\email{jakob.stoye@tu-braunschweig.de})
}
}

\usepackage{amsopn}

%% file: InjEuclideanStiefel.bbl
\begin{thebibliography}{10}
	
	\bibitem{Absilbook:2008}
	{\sc P.-A. Absil, R.~Mahony, and R.~Sepulchre}, {\em Optimization Algorithms on
		Matrix Manifolds}, Princeton University Press, Princeton, NJ, 2008.
	
	\bibitem{absil2024ultimate}
	{\sc P.~A. Absil and S.~Mataigne}, {\em The ultimate upper bound on the
		injectivity radius of the {S}tiefel manifold}, 2024,
	\url{https://arxiv.org/abs/2403.02079}.
	
	\bibitem{BennerGugercinWillcox2015}
	{\sc P.~Benner, S.~Gugercin, and K.~Willcox}, {\em A survey of projection-based
		model reduction methods for parametric dynamical systems}, SIAM Review, 57
	(2015), pp.~483--531, \url{https://doi.org/10.1137/130932715},
	\url{http://dx.doi.org/10.1137/130932715},
	\url{https://arxiv.org/abs/http://dx.doi.org/10.1137/130932715}.
	
	\bibitem{boumal2023}
	{\sc N.~Boumal}, {\em An Introduction to Optimization on Smooth Manifolds},
	Cambridge University Press, Cambridge, 2023.
	
	\bibitem{Celledoni_2020}
	{\sc E.~Celledoni, S.~Eidnes, B.~Owren, and T.~Ringholm}, {\em Energy
		preserving methods on {R}iemannian manifolds}, Mathematics of Computation,
	(2020), pp.~699--716.
	
	\bibitem{Chakraborty2018}
	{\sc R.~Chakraborty and B.~Vemuri}, {\em Statistics on the (compact) {S}tiefel
		manifold: Theory and applications}, The Annals of Statistics, 47 (2018),
	\url{https://doi.org/10.1214/18-AOS1692}.
	
	\bibitem{DoCarmo2013riemannian}
	{\sc M.~P. do~Carmo}, {\em {R}iemannian Geometry}, Mathematics: Theory \&
	Applications, Birkh{\"a}user Boston, 1992,
	\url{https://books.google.de/books?id=ct91XCWkWEUC}.
	
	\bibitem{EdelmanAriasSmith:1999}
	{\sc A.~Edelman, T.~A. Arias, and S.~T. Smith}, {\em The geometry of algorithms
		with orthogonality constraints}, SIAM Journal on Matrix Analysis and
	Applications, 20 (1998), pp.~303--353,
	\url{https://doi.org/10.1137/S0895479895290954},
	\url{http://dx.doi.org/10.1137/S0895479895290954}.
	
	\bibitem{gallier2011geometric}
	{\sc J.~Gallier}, {\em Geometric Methods and Applications: For Computer Science
		and Engineering}, Texts in Applied Mathematics, Springer New York, 2011,
	\url{https://doi.org/https://doi.org/10.1007/978-1-4613-0137-0}.
	
	\bibitem{Gluck:1966}
	{\sc H.~Gluck}, {\em Higher curvatures of curves in {E}uclidean space}, The
	American Mathematical Monthly, 73 (1966), pp.~699--704.
	
	\bibitem{HueperRollingStiefel2008}
	{\sc K.~H\"uper, M.~Kleinsteuber, and F.~Silva~Leite}, {\em Rolling {S}tiefel
		manifolds}, International Journal of Systems Science, 39 (2008),
	pp.~881--887, \url{https://doi.org/10.1080/00207720802184717},
	\url{https://arxiv.org/abs/https://doi.org/10.1080/00207720802184717}.
	
	\bibitem{HueperMarkinaLeite2020}
	{\sc K.~H{\"u}per, I.~Markina, and F.~Silva~Leite}, {\em A {L}agrangian
		approach to extremal curves on {S}tiefel manifolds}, Journal of Geometrical
	Mechanics, 13 (2021), pp.~55--72.
	
	\bibitem{Klingenberg:1978}
	{\sc W.~Klingenberg}, {\em A Course in Differential Geometry}, no.~51 in
	Graduate Texts in Mathematics, Springer-Verlag, New York, 1978.
	
	\bibitem{kuhnel:2015}
	{\sc W.~K{\"u}hnel}, {\em Differential Geometry: Curves -- Surfaces --
		Manifolds}, Student Mathematical Library, American Mathematical Society,
	3~ed., 2015, \url{https://books.google.dk/books?id=qNBYCwAAQBAJ}.
	
	\bibitem{Lee:RiemannianManifolds}
	{\sc J.~M. Lee}, {\em Riemannian Manifolds: An Introduction to Curvature},
	Graduate Texts in Mathematics, Springer New York, NY, 1997,
	\url{https://doi.org/https://doi.org/10.1007/b98852}.
	
	\bibitem{Monterde:2005}
	{\sc J.~Monterde}, {\em Curves with constant curvature ratios}, Boletín de la
	Sociedad Matemática Mexicana, 13 (2005).
	
	\bibitem{nguyen2022curvature}
	{\sc D.~Nguyen}, {\em Curvatures of {S}tiefel manifolds with deformation
		metrics}, Journal of Lie Theory, 32 (2022), pp.~563--600.
	
	\bibitem{FLETCHER2020}
	{\sc X.~Pennec, S.~Sommer, and T.~Fletcher}, eds., {\em {R}iemannian Geometric
		Statistics in Medical Image Analysis}, Academic Press, 2020.
	
	\bibitem{petersen2016riemannian}
	{\sc P.~Petersen}, {\em Riemannian Geometry}, Graduate Texts in Mathematics,
	Springer International Publishing, 2016,
	\url{https://books.google.dk/books?id=YwbNCwAAQBAJ}.
	
	\bibitem{sato2021}
	{\sc H.~Sato}, {\em {R}iemannian Optimization and Its Applications},
	SpringerBriefs in Electrical and Computer Engineering, Springer International
	Publishing, 2021, \url{https://books.google.dk/books?id=f9oeEAAAQBAJ}.
	
	\bibitem{Turaga_2008}
	{\sc P.~K. Turaga, V.~A., and R.~Chellappa}, {\em Statistical analysis on
		{S}tiefel and {G}rassmann manifolds with applications in computer vision}, in
	2008 IEEE Conference on Computer Vision and Pattern Recognition, June 2008,
	pp.~1--8, \url{https://doi.org/10.1109/CVPR.2008.4587733}.
	
	\bibitem{Wong1967}
	{\sc Y.-C. Wong}, {\em Differential geometry of {G}rassmann manifolds},
	Proceedings of the National Academy of Sciences of the United States of
	America, 57 (1967), pp.~589--594.
	
	\bibitem{Wong1968}
	{\sc Y.-C. Wong}, {\em Sectional curvatures of {G}rassmann manifolds},
	Proceedings of the National Academy of Sciences of the United States of
	America, 60 (1968), pp.~75--79, \url{http://www.jstor.org/stable/58432}.
	
	\bibitem{Zim21}
	{\sc R.~Zimmermann}, {\em Manifold interpolation}, in Volume 1 System- and
	Data-Driven Methods and Algorithms, P.~Benner, S.~Grivet-Talocia,
	A.~Quarteroni, G.~Rozza, W.~Schilders, and L.~M. Silveira, eds., De Gruyter,
	2021, pp.~229--274, \url{https://doi.org/10.1515/9783110498967-007}.
	
	\bibitem{ZimmermannHueper2022}
	{\sc R.~Zimmermann and K.~H\"{u}per}, {\em Computing the {R}iemannian logarithm
		on the {S}tiefel manifold: Metrics, methods, and performance}, SIAM Journal
	on Matrix Analysis and Applications, 43 (2022), pp.~953--980,
	\url{https://doi.org/10.1137/21M1425426}.
	
	\bibitem{zimmermannstoye:2024}
	{\sc R.~Zimmermann and J.~Stoye}, {\em High curvature means low rank: On the
		sectional curvature of {G}rassmann and {S}tiefel manifolds and the underlying
		matrix trace inequalities}, 2024, \url{https://arxiv.org/abs/2403.01879}.
	
\end{thebibliography}
